\PassOptionsToPackage{dvipsnames}{xcolor}
\documentclass[12pt,onecolumn]{IEEEtran} 

\usepackage{amsmath} 
\usepackage{amssymb}  
\usepackage{amsbsy}
\usepackage{xcolor}
\usepackage[T1]{fontenc}
\usepackage[latin1]{inputenc}
\usepackage{varioref}
\usepackage{subfigure}
\usepackage{eurosym}
\usepackage{cite}
\usepackage{url}
\usepackage{dblfloatfix}
\usepackage{booktabs,fixltx2e}
\usepackage[flushleft]{threeparttable}

\usepackage{natbib}            
\usepackage{graphicx}          

\usepackage{acronym}

\usepackage[detect-all]{siunitx}

\title{Impact of Low Rotational Inertia on\\Power System Stability and Operation}

\author{
	\IEEEauthorblockN{Andreas Ulbig, Theodor S. Borsche and Göran Andersson\\}  
	\IEEEauthorblockA{Power Systems Laboratory, ETH Zurich\\ \textnormal{ulbig}$\;\textbar\;$\textnormal{borsche}$\;\textbar\;$\textnormal{andersson}$\ $\textnormal{@$\,$eeh.ee.ethz.ch}}
}

\acrodef{phs}[PHS]{pumped hydro storage}
\acrodef{ror}[ROR]{run-of-river hydro}
\acrodef{hdd}[HDD]{heating degree day}
\acrodef{cdd}[CDD]{cooling degree day}
\acrodef{vpp}[VPP]{virtual power plants}
\acrodef{epex}[EPEX]{European Power Exchange}
\acrodef{eex}[EEX]{European Energy Exchange}
\acrodef{res}[RES]{Renewable Energy Sources}
\acrodef{fit}[FIT]{Feed-In Tariffs}
\acrodef{eeg}[EEG]{"Erneuerbare-Energien-Gesetz"}
\acrodef{pv}[PV]{photovoltaics}
\acrodef{tso}[TSO]{Transmission System Operator}
\acrodef{kti}[KTI]{Swiss Innovation Promotion Agency}
\acrodef{mape}[MAPE]{mean absolute prediction error}
\acrodef{mse}[MSE]{mean square error}
\acrodef{dlr}[DLR]{dynamic line rating}
\acrodef{dr}[DR]{Demand Response}
\acrodef{dsm}[DSM]{Demand-Side Management}
\acrodef{dsp}[DSP]{Demand-Side Participation}
\acrodef{rps}[RPS]{Renewable Portfolio Standards}
\acrodef{chp}[CHP]{Combined Heat and Power}
\acrodef{soc}[SOC]{State-of-Charge}
\acrodef{nreap}[NREAP]{National Renewable Energy Action Plan}
\acrodef{csp}[CSP]{Concentrating Solar Thermal Power}

\usepackage{hyperref}
\usepackage{bookmark}

\begin{document}
\definecolorset{cmyk}{}{}{%
    eth1, 1,    0.7,  0,   0.30;%
    eth2, 0.75, 0.4,  1,   0.40;%
    eth3, 1,    0.5,  0,   0;%
    eth4, 0.3,  0,    1,   0.55;%
    eth5, 0.2,  1,    0,   0.20;%
    eth6, 0,    0,    0,   0.7;%
    eth7, 0,    0.9,  0.8, 0.20;%
    eth8, 1,    0.25, 0.3, 0.1;%
    eth9, 0,    0.55, 1,   0.40;%
    eth10,0.6,  0,    1,   0 }

\definecolorset{RGB}{}{}{%
    eth1rgb,  31,  64, 122;%
    eth2rgb,  72,  90,  44;%
    eth3rgb,  18, 105, 176;%
    eth4rgb, 114, 121,  28;%
    eth5rgb, 145,   5, 106;%
    eth6rgb, 111, 111, 100;%
    eth7rgb, 168,  50,  45;%
    eth8rgb,   0, 122, 150;%
    eth9rgb, 149,  96,  19;%
    eth10rgb,140, 182,  60}
\maketitle

\IEEEpeerreviewmaketitle   					                      	

\begin{abstract}                          
Large-scale deployment of \ac{res} has led to significant generation shares of variable \ac{res} in power systems worldwide. 
\ac{res} units, notably inverter-connected wind turbines and \ac{pv} that as such do not provide rotational inertia, are effectively displacing conventional generators and their rotating machinery. The traditional assumption that grid inertia is sufficiently high with only small variations over time is thus not valid for power systems with high RES shares. This has implications for frequency dynamics and power system stability and operation. Frequency dynamics are faster in power systems with low rotational inertia, making frequency control and power system operation more challenging. 

This paper investigates the impact of low rotational inertia on power system stability and operation, contributes new analysis insights and offers mitigation options for low inertia impacts.
\end{abstract}

\vspace{0.25cm}


%


\section{Introduction}
Traditionally, power system operation is based on the assumption that electricity generation, in the form of thermal power plants, reliably supplied with fossil or nuclear fuels, or hydro plants, is fully dispatchable, i.e.~controllable, and involves rotating synchronous generators. Via their stored kinetic energy they add rotational inertia, an important property of frequency dynamics and stability. The contribution of inertia is an inherent and crucial feature of rotating synchronous generators. Due to electro-mechanical coupling, a generator's rotating mass provides kinetic energy to the grid (or absorbs it from the grid) in case of a frequency deviation~$\Delta f$. The kinetic energy provided is proportional to the rate of change of frequency~$\Delta \dot{f}$~\cite{Kundur1994}. 
The grid frequency $f$ is directly coupled to the rotational speed of a synchronous generator and thus to the active power balance. 
Rotational inertia, i.e.~the inertia constant $H$, minimizes~$\Delta \dot{f}$ in case of frequency deviations. This renders frequency dynamics more benign, i.e~slower, and thus increases the available response time to react to fault events such as line losses, power plant outages or large-scale set-point changes of either generation or load units.

Maintaining the grid frequency within an acceptable range is a necessary requirement for the stable operation of power systems. Frequency stability and in turn also stable operation both depend on the active power balance, meaning that the total power feed-in minus the total load consumption (including system losses) is kept close to zero. In normal operation small variations of this balance occur spontaneously. Deviations from its nominal value $f_0$, e.g.~50~Hz or 60~Hz depending on region, should be kept small, as damaging vibrations in synchronous machines and load shedding occur for larger deviations. 
This can influence the whole power system, in the worst case ending in fault cascades and black-outs.
Low levels of rotational inertia in a power system, caused in particular by high shares of inverter-connected \ac{res}, i.e.~wind turbine and \ac{pv} units that normally do not provide any rotational inertia, have implications on frequency dynamics. They are becoming faster in power systems with low rotational inertia. This can lead to situations in which traditional frequency control schemes become too slow with respect to the disturbance dynamics for preventing large frequency deviations and the resulting consequences. 
The loss of rotational inertia and its increasing time-variance lead to new frequency instability phenomena in power systems. Frequency and power system stability may be at risk.

An exemplary analysis of the German power system shows the relevance of the above mentioned trends. Throughout the year 2012 there have been several occasions hours in which around 50\% of overall load demand was covered by wind\&\ac{pv} units. 
The regional inertia within the German power system dropped to significantly lower levels than usual due to the temporary lack of dispatched conventional generators and their rotating machinery. 
With the increase of inverter-connected RES generation, low inertia situations will become more widespread and with it faster frequency dynamics and the associated operational risks.  

\begin{figure*}[b]
  \centering               
   {\includegraphics[trim = 1.0cm   0.5cm   1.5cm   1.0cm, clip=true, angle=0, width=0.46\linewidth, keepaspectratio, draft=false]{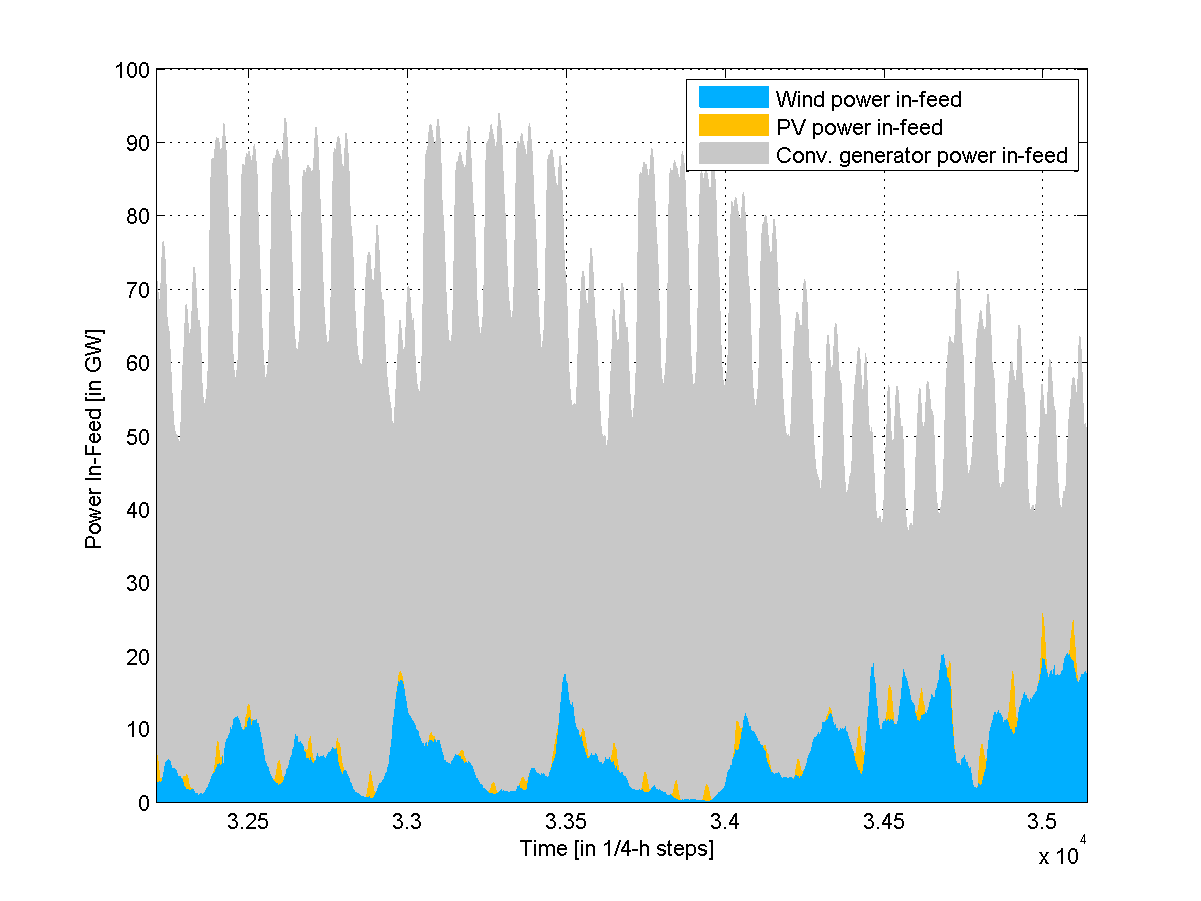}}
  \centering               
   {\includegraphics[trim = 1.0cm   0.5cm   1.5cm   1.0cm, clip=true, angle=0, width=0.46\linewidth, keepaspectratio, draft=false]{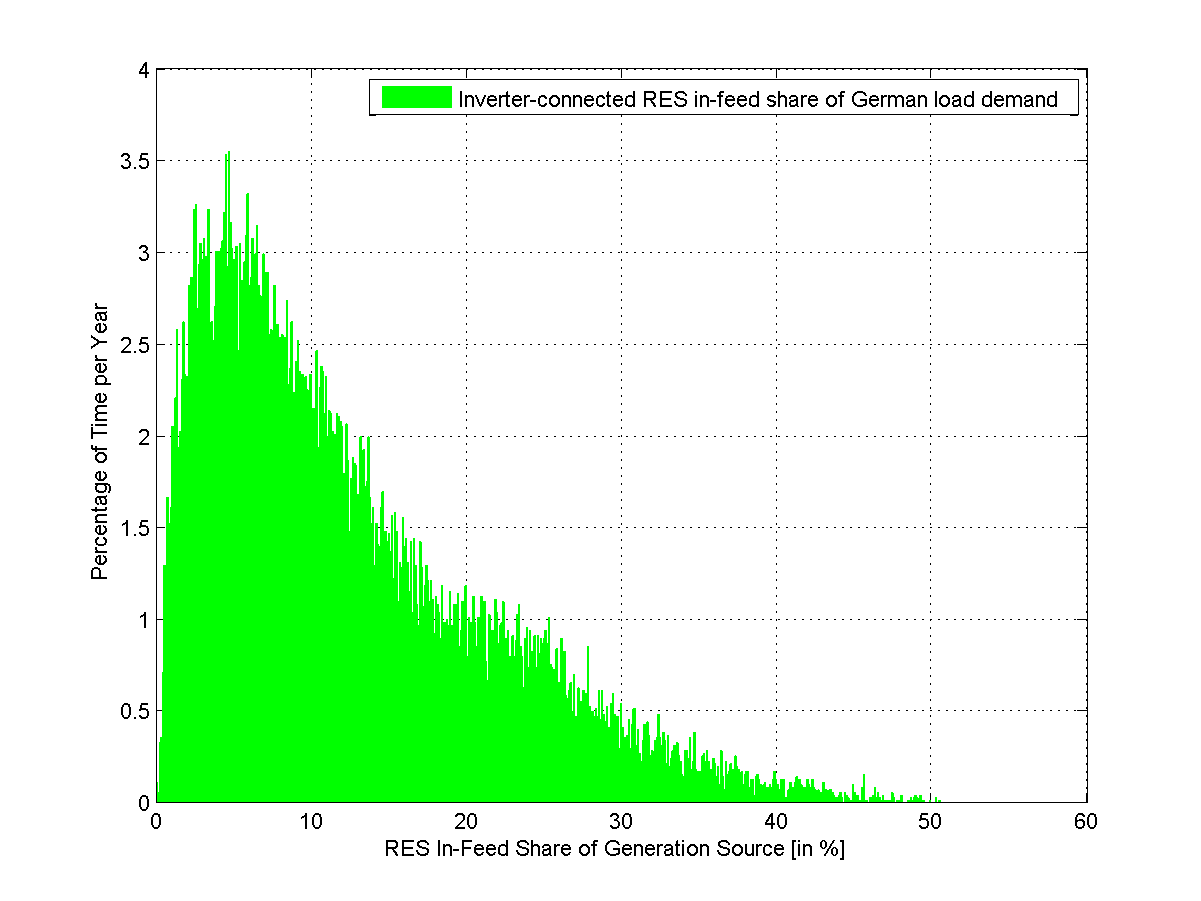}}
  \caption{(a)~Power Dispatch Situation in German Power System (December~2012). (b)~Histogram of Inverter-Connected Power Feed-in Shares in German Power System (full-year~2012).} 
  \label{fig:RESShare}
\end{figure*}

The remainder of this paper is organized as follows: Section~\ref{sec:RES} discusses the rapid large-scale deployment of RES generation in many countries and the arising challenges for power system operation. Section~\ref{sec:Inertia} explains rotational inertia in more detail and assesses to what extent inverter-connected generation units reduce inertia and render it time-variant. This is followed by an analysis of the impacts of reduced inertia on power system stability in Section~\ref{sec:Stability} and power system operation in Section~\ref{sec:Operation}. Finally, a conclusion and an outlook are given in Section~\ref{sec:Conclusion}.

\section{Impacts of Rising Renewable Energy Shares for Power System Operation}\label{sec:RES}

Facing the challenge of having to reduce $\mathrm{CO_2}$ emissions due to climate change concerns as well as security of supply issues of fossil fuels, many countries nowadays are committed to increasing the share of renewable energy sources (RES) in their electric power systems.

Large-scale deployment of \ac{res} generation, notably in the form of wind turbines and \ac{pv} units, ranging from small and highly distributed units, e.g.~roof-top PV with a rating of a few kilowatts (kW), to large units, e.g.~large PV and wind farms with hundreds of megawatts (MW), has led to significant generation shares of variable \ac{res} power feed-in in power systems worldwide.  
\ac{res} capacity comprised about 25\% of total global power generation capacity and produced an estimated 20.3\% of global electricity demand by by year-end 2011. Although most \ac{res} electricity is still provided by hydro power ($15\, \%$) other renewables ($5.3\,\%$) are on the rise. Of the world's total generation capacity estimated at $5360\,\textrm{GW}_\textrm{el}$ by year-end~2011, wind power made up $238\,\textrm{GW}_\textrm{el}$ ($4.4\,\%$), solar \ac{pv} $70\,\textrm{GW}_\textrm{el}$ ($1.3\,\%$) whereas \ac{csp} only contributed $1.8\,\textrm{GW}_\textrm{el}$ ($0.03\,\%$). In the European Union (EU-28), with a total generation capacity of around $870\,\textrm{GW}_\textrm{el}$, wind power made up $94\,\textrm{GW}_\textrm{el}$ ($10.8\,\%$) and solar \ac{pv} $51\,\textrm{GW}_\textrm{el}$ ($5.9\,\%$)~\cite{REN21:2012}.

In Germany, the \ac{res} share of electricity generation increased from 4.7\% of net load demand in~1998 to more than 20\% in~2012. \ac{res} generation capacity is dominated by wind, \ac{pv} and hydro generation with an absolute share of net load demand of 8.3\%, 5.0\% and 3.9\%, respectively, in 2012. The remainder was made up of biomass, land-fill and bio gas generation (3--4\%)~\cite{RESfigure}.

Due to the rising \ac{res} shares, the number of hours per year in which \ac{res} feed-in makes up a large part or the majority of power production in a grid region is also increasing.
This is illustrated for the case of Germany in Fig.~\ref{fig:RESShare}~(a)--(b). There, the power dispatch situation of wind\&\ac{pv} units and conventional generation in the German power system is illustrated for December~2012 (31~days). Also, the histogram of the total inverter-connected \ac{res} feed-in, i.e.~wind\&\ac{pv}, as a share of the total load demand in Germany is given for the full year~2012. In this particular year the share of inverter-connected \ac{res} units often reached significant levels: a share of 30\% or more was reached for 495~hours a year (5.6\%), 40\% or more for 221~hours (2.5\%) and a record 50\% for 0.75~hours (0.009\%), respectively. 


\section{Time-Variance of Grid Inertia}\label{sec:Inertia}

\begin{figure*}[b]
  \centering               
   {\includegraphics[trim = 1.0cm   0.5cm   1.5cm   1.0cm, clip=true, angle=0, width=0.46\linewidth, keepaspectratio, draft=false]{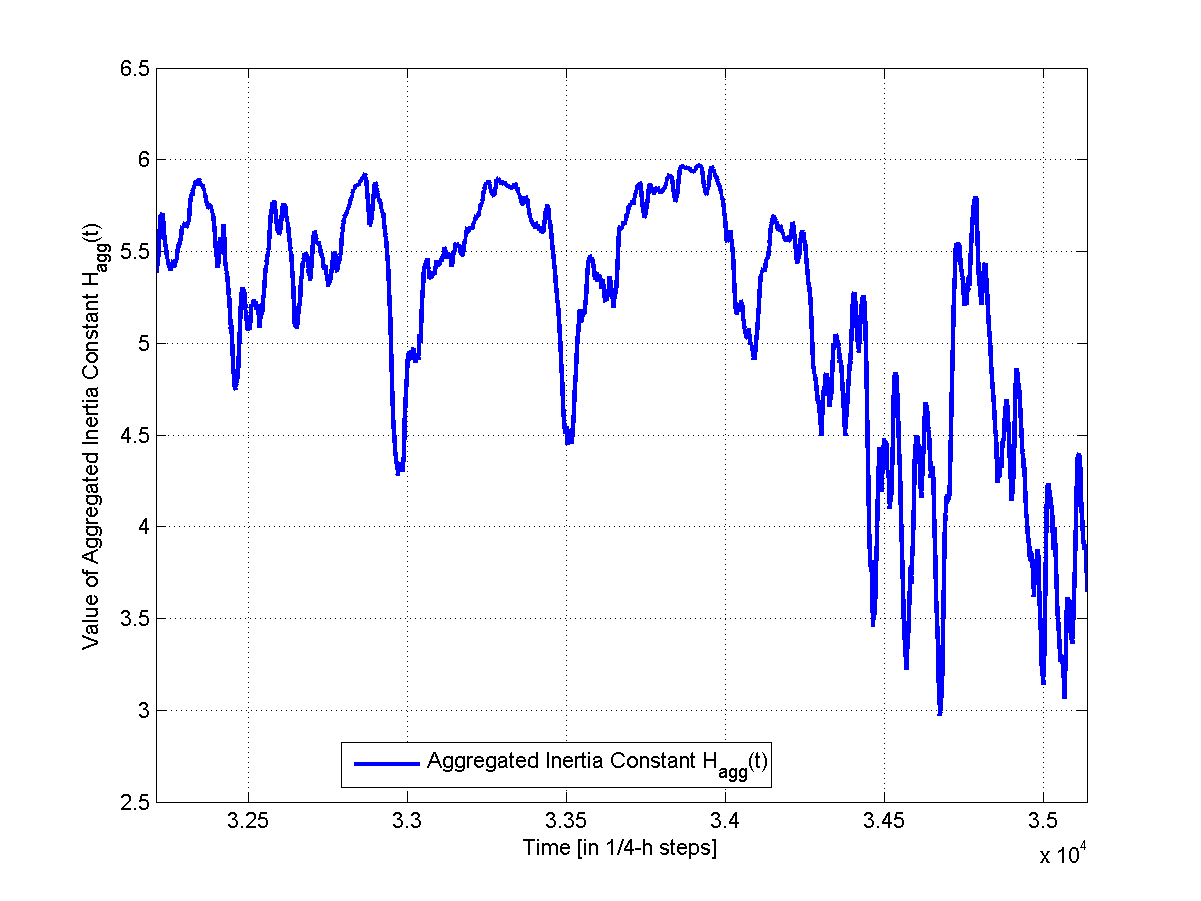}}
   {\includegraphics[trim = 1.0cm   0.5cm   1.5cm   1.0cm, clip=true, angle=0, width=0.46\linewidth, keepaspectratio, draft=false]{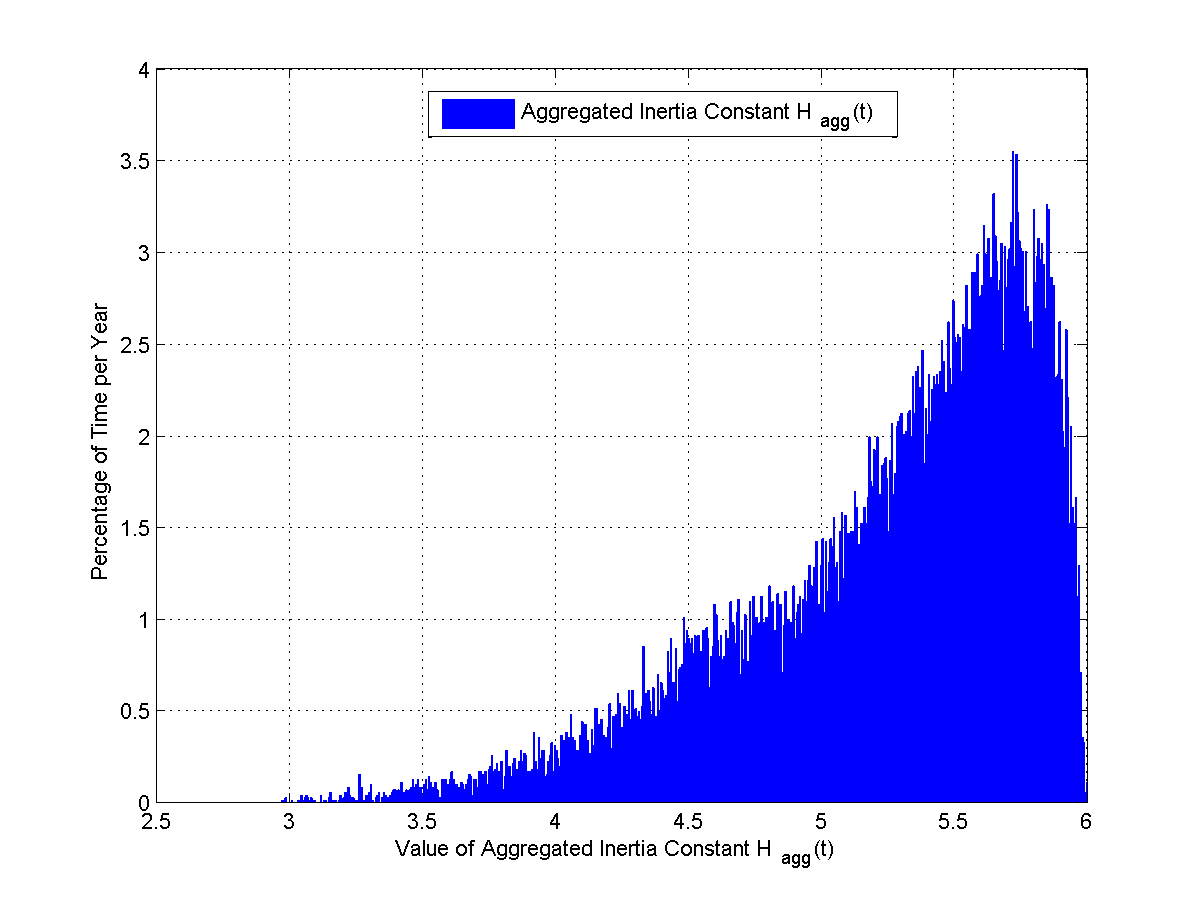}}
  \caption{(a)~Time-Variant Aggregated Rotational Inertia $H_{\textrm{agg}}$ in German Power System (December~2012). It is assumed that conventional generators provide inertia ($H_{\textrm{conv}} = 6\,\textrm{s}$) and inverter-connected RES generators do not ($H_{\textrm{RES}} = 0\,\textrm{s}$). (b)~Histogram of Aggregated Rotational Inertia in German Power System (full-year~2012).}
  \label{fig:Inertia} 
\end{figure*}

In the following the basic modeling concepts for rotational inertia in power systems as well as synchronous power systems in general are presented.
 
\subsection{Modeling Inertial Response}

Following a frequency deviation, kinetic energy stored in the rotating masses of the generator system is released, rendering power system frequency dynamics slower and, hence, easier to regulate. The rotational energy is given as
\begin{equation}
E_{\textrm{kin}} = \frac{1}{2}J(2\pi f_\textrm{m})^{2} \ \ ,\label{eq:kinENERGY}\end{equation}
with $J$ as the moment of inertia of the synchronous machine and $f_{\textrm{m}}$ the rotating frequency of the machine. The inertia constant $H$ for a synchronous machine is defined by
\begin{equation}
H = \frac{E_{\textrm{kin}}}{S_\textrm{B}} = \frac{J (2\pi f_{\textrm{m}})^{2}}{2S_\textrm{B}} \ \ ,
\end{equation}
with $S_\textrm{B}$ as the rated power of the generator and $H$ denoting the time duration during which the machine can supply its rated power solely with its stored kinetic energy. Typical values for $H$ are in the range of $2$--$10$~s \cite[Table 3.2]{Kundur1994}. The classical swing equation, a well-known model representation for synchronous generators, describes the inertial response of the synchronous generator as the change in rotational frequency $f_\textrm{m}$ (or rotational speed $\omega_\textrm{m} = 2\pi \cdot f_\textrm{m}$) of the synchronous generator following a power imbalance as
\begin{equation}
\dot{E}_{\textrm{kin}} = J (2 \pi)^{2} f_{\textrm{m}} \cdot \dot{f}_\textrm{m}  = \frac{2HS_\textrm{B}}{f_{\textrm{m}}} \cdot \dot{f}_{\textrm{m}} = (P_{\textrm{m}} - P_{\textrm{e}}) \ \ , \label{eq:swing equation}
\end{equation}
with $P_{\textrm{m}}$ as the mechanical power supplied by the generator and $P_{\textrm{e}}$ as the electric power demand.

Noting that frequency excursions are usually small deviations around the reference value, we replace $f_\textrm{m}$ by $f_0$ and $P_\textrm{m}$ by $P_\textrm{m\,0}$, and complete the classical Swing Equation by adding frequency-dependent load damping, a self-stabilizing property of power systems, by formulating
\begin{equation}
\dot{f}_\textrm{m} = -\frac{f_{0}}{2 H S_\textrm{B} D_{\textrm{load}}}f_\textrm{m} + \frac{f_{0}}{2HS_\textrm{B}}(P_\textrm{m,\,0}-P_\textrm{e}) \ \ .\label{eq:swing control}
\end{equation}

Here $f_{0}$ is the reference frequency and $D_\textrm{load}$ denotes the frequency-dependent load damping constant. $P_\textrm{m,\,0}$ is the nominally scheduled mechanical generator power. Another definition of load damping is $k_\textrm{load}$ with $k_\textrm{load} = \frac{1}{D_\textrm{load}}$. Please note that in literature, concurrent labelings like $D_\textrm{l}$ (or simply $D$) and $k_\textrm{l}$ (or $k$) are also in wide use.
The high share of conventional generators is translated into a large rotational inertia of the here presented power system. The higher the inertia constant $H$, the slower and more benign are frequency dynamics, i.e.~for identical faults frequency deviations $f_m$ and their derivatives $\dot{f}_m$ are smaller. 

With an increasing penetration of inverter-connected power units, the rotational inertia of power systems is reduced and becomes highly time-variant as wind\&\ac{pv} shares are fluctuating heavily throughout the year. This is notably a concern for small power networks, e.g.~island or micro grids, with a high share of generation capacity not contributing any inertia as was discussed and illustrated, for example, in~\cite{VanHertem2012}. Frequency stabilization becomes thus more difficult. Appropriate adaptations of grid codes are hence needed.

\subsection{Aggregated Swing Equation Model}

Modeling interconnected power systems, i.e.~different aggregated generator and load nodes that are connected via tie-lines, can be realized in a similar fashion as modeling individual generators. Reformulating the classical Swing Equation (Eq.~\ref{eq:swing control}) for a power system with $n$~generators, $j$~loads and $l$~connecting tie-lines, leads to the so-called Aggregated Swing Equation~(ASE)~\cite{Kundur1994}
\begin{equation}\label{eq:ASE}
    \dot{f} = - \frac{f_{0}}{2 H S_{\textrm{B}} D_{\textrm{load}}} f + \frac{f_{0}}{2HS_\textrm{B}} ( P_\textrm{m} - P_\textrm{load} - P_\textrm{loss} )\ ,
\end{equation}
with
\begin{eqnarray*}
f  = \frac{\sum_{i=1}^n H_{i} \, S_\textrm{B,i} \, f_{i}}{\sum_{i=1}^n H_{i} \, S_\textrm{B,i}} \, , \
S_\textrm{B}  = \sum_{i=1}^n S_\textrm{B,i} \, , \
H  = \frac{\sum_{i=1}^n H_{i} S_\textrm{B,i} } {S_\textrm{B}} \, , \\
P_\textrm{m} = \sum_{i=1}^n P_{\textrm{m},i} \, , \
P_\textrm{load} = \sum_{i=1}^j P_{\textrm{load},i} \, , \
P_\textrm{loss} = \sum_{i=1}^l P_{\textrm{loss},i} \, .
\end{eqnarray*}

Here the term $f$ is the Center of Inertia (COI) grid frequency, $H$~the aggregated inertia constant of the $n$ generators, $S_\textrm{B}$~the total rated power of the generators, $P_\textrm{m}$ the total mechanical power of the generators, $P_\textrm{load}$ the total system load of the grid and $P_\textrm{loss}$ the total transmission losses of the $l$ lines making up the grid topology and $f_{0}=50$~Hz. The term $D_\textrm{load}$ is the frequency damping of the system load, which is assumed here to be constant and uniform. 
All power system parameters are given in Table~\ref{tab:parameters}. 

The ASE model (Eq. \ref{eq:ASE}) is valid for a highly meshed grid, in which all units can be assumed to be connected to the same grid bus, representing the Center of Inertia of the given grid. Since load-frequency disturbances are normally relatively small, linearized swing equations with $\Delta f_{i}=f_{i}-f_{0}$ can be used. Considering the system change ($\Delta$) before and after a disturbance, the relative formulation of the ASE, assuming that $\Delta P_\textrm{loss} = 0$, is
\begin{equation}
\label{eq:DeltaASE}
\Delta\dot{f} = -\frac{f_{0}}{2 H S_\textrm{B} D_\textrm{load}}\Delta f + \frac{f_{0}}{2HS_\textrm{B}} \left( \Delta P_\textrm{m} - \Delta P_\textrm{load} \right) \quad.\hspace{-0.1cm}
\end{equation}

In frequency stability analysis often the assumption is used that the (aggregated) inertia constant $H$ is constant (and the same) for all swing equations of a multi-area system. This assumption was valid in the past but is nowadays increasingly tested by reality as is illustrated in Fig.~\ref{fig:Inertia}, again for the case of the German power system. It shows that its aggregated inertia $H_{\textrm{agg}}$, as calculated using the respective equation in~(\ref{eq:ASE}), has indeed become highly time-variant and fluctuates between its nominal value of 6~s, i.e.~at times when only conventional generators are dispatched, and significantly lower levels of 3--4~s, i.e.~at times when significant shares of wind\&\ac{pv} generation are deployed. The lowest level of rotational inertia of this year was reached during the Christmas vacation in which demand levels were at their lowest (in December~2012), while notably wind power feed-in was unusually high.
The histogram for the full year 2012 reveals that inertia levels drop to rather low levels for a significant part of the time: $H_{\textrm{agg}}$ was below 4~s for 293~hours (3.3\%) and below 3.5~s for 57~h (0.65\%) of the time.
The qualitative results of this example are valid also for the inertia situation in other countries with high \ac{res} shares.

As this section and the previous one show, coping with the fluctuating electricity production from variable RES,~i.e. wind turbines and PV, is a challenge for the operation of electric power systems in many aspects. The increasing share of inverter-based power generation and the associated displacement of usually large-scale and fully controllable generation units and their rotational masses, in particular has the following consequences:
\begin{enumerate}
\item{The pool of suitable conventional power plants for providing traditional control reserve power is significantly diminished.}

\item{The rotational inertia of power systems becomes markedly time-variant and is reduced, often non-uniformly within the grid topology, as will be presented in the following section.}
\end{enumerate}

\section{Impact of Low Rotational Inertia on Power System Stability}\label{sec:Stability}

Frequency dynamics of single-area as well as multi-area power systems are usually modeled and analyzed employing the Swing Equation approach introduced in Section~\ref{sec:Inertia}.

It is known that frequency dynamics for a system with $n$~areas can become chaotic in case~$n \ge 3$; confer to~\cite{Kopell1982},~\cite{Varaiya1987} or~\cite{Berggren1993} for more details. Analyzing the stability properties of swing equation models of power systems constituted a sizeable research stream in the~1980s and early~1990s. Although the analysis presented back then assumed that rotational inertia constants could vary from one grid region to another, its time-variance caused by massive inverter-connected \ac{res} feed-in was not considered at the time as only very few wind\&\ac{pv} units existed.

The following analyses use a three-area power system that was simplified to a two-area model, as the reference voltage angle and frequency of the third grid area are kept at zero, i.e.~$\delta_3=0,\ \omega_3=0$. The modeling is based on the Swing Equation approach and follows the line of thought presented in the work of~\cite{Varaiya1987}: 
\begin{eqnarray}\label{eq:Phase}
 \dot{\delta}_1  &=& \omega_1 \\
 \dot{\delta}_2  &=& \omega_2 \nonumber \\
 \dot{\omega}_1 &=&  \frac{1}{M_1} \left[ \Delta P_1 - k_1 \omega_1 - V_1 V_2 B_2 \sin(\delta_1 - \delta_2) \right. \nonumber \\
 &&- \left. V_1 V_3 B_3 \sin(\delta_1-\delta_3) \right] \nonumber \\
 \dot{\omega}_2 &=&  \frac{1}{M_2} \left[ \Delta P_2 - k_2 \omega_2 - V_2 V_1 B_1 \sin(\delta_2 - \delta_1) \right. \nonumber \\
 &&- \left. V_2 V_3 B_3 \sin(\delta_2-\delta_3) \right] \nonumber \,. 
\end{eqnarray}

Here the voltage levels $V_i$ are assumed to be nominal, i.e.~1~p.u. The specifications of all parameters are given in Table~\ref{tab:parameters}. Also, the familiar terms for rotational inertia and power deviations are linked with the previous equations introduced in Section~\ref{sec:Inertia} via 
\begin{equation}
M_i = \frac{2 H_i S_{\textrm{B}_i}}{2 \pi f_0} = J_i \omega_i \ \textrm{and} \ \Delta P_i = \left( \Delta P_\textrm{m} - \Delta P_\textrm{load} \right) \quad.
\end{equation}
For the sake of simplicity in presentation and for easier comparison with the related work previously mentioned, we use in this section the inertia constants $M_i$ instead of $H_i$ and the angular frequency $\omega_i = 2\pi f_i$, given in rad/s, instead of using directly the frequency $f_i$.
From Eq.~\ref{eq:Phase} one can analytically deduce that the inertia constant $M_i$ mitigates the impact of \emph{shocks} such as sudden power faults $\Delta P_i$ on the angular frequency $\omega_i$ since $\dot{\omega}_i \sim \left( \frac{\Delta P_i}{M_i} \right)$. Also the damping coefficient $k_i$ has a stabilizing effect on $\omega_i$ since $\dot{\omega}_i \sim \left(-\frac{k_i}{M_i} \cdot \omega_i \right)$. The achieved stabilizing effect depends, however, also on the ratio of $\left(\frac{k_i}{M_i}\right)$. 
Both the values of the inertia constant $M_i$ and the damping coefficient $k_i$ are thus vital for power system stability.

\begin{figure}[b!]
  \centering               
{\includegraphics[trim = 1.0cm   0cm   1.5cm   0.0cm, clip=true, angle=0, width=0.36\linewidth, keepaspectratio, draft=false]{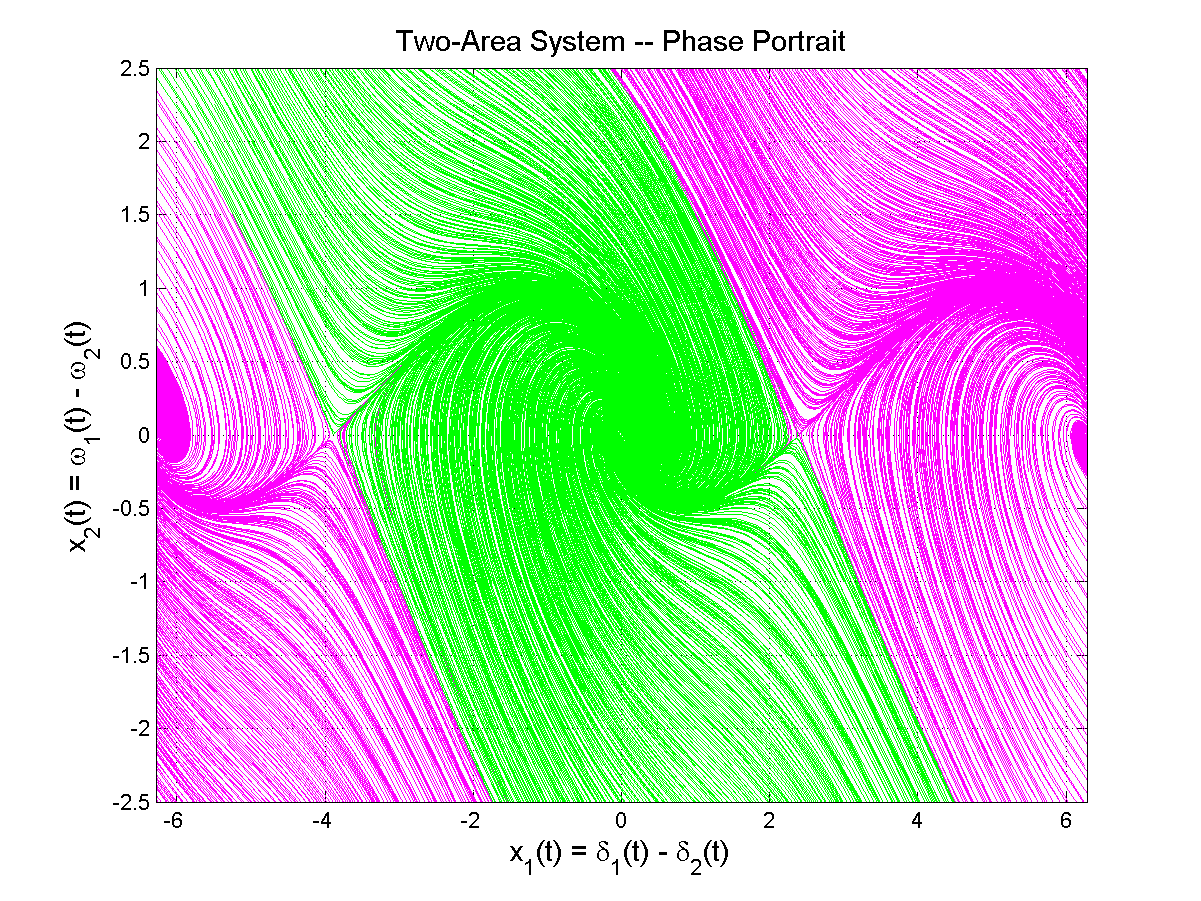}}
{\includegraphics[trim = 1.0cm   0cm   1.5cm   0.0cm, clip=true, angle=0, width=0.36\linewidth, keepaspectratio, draft=false]{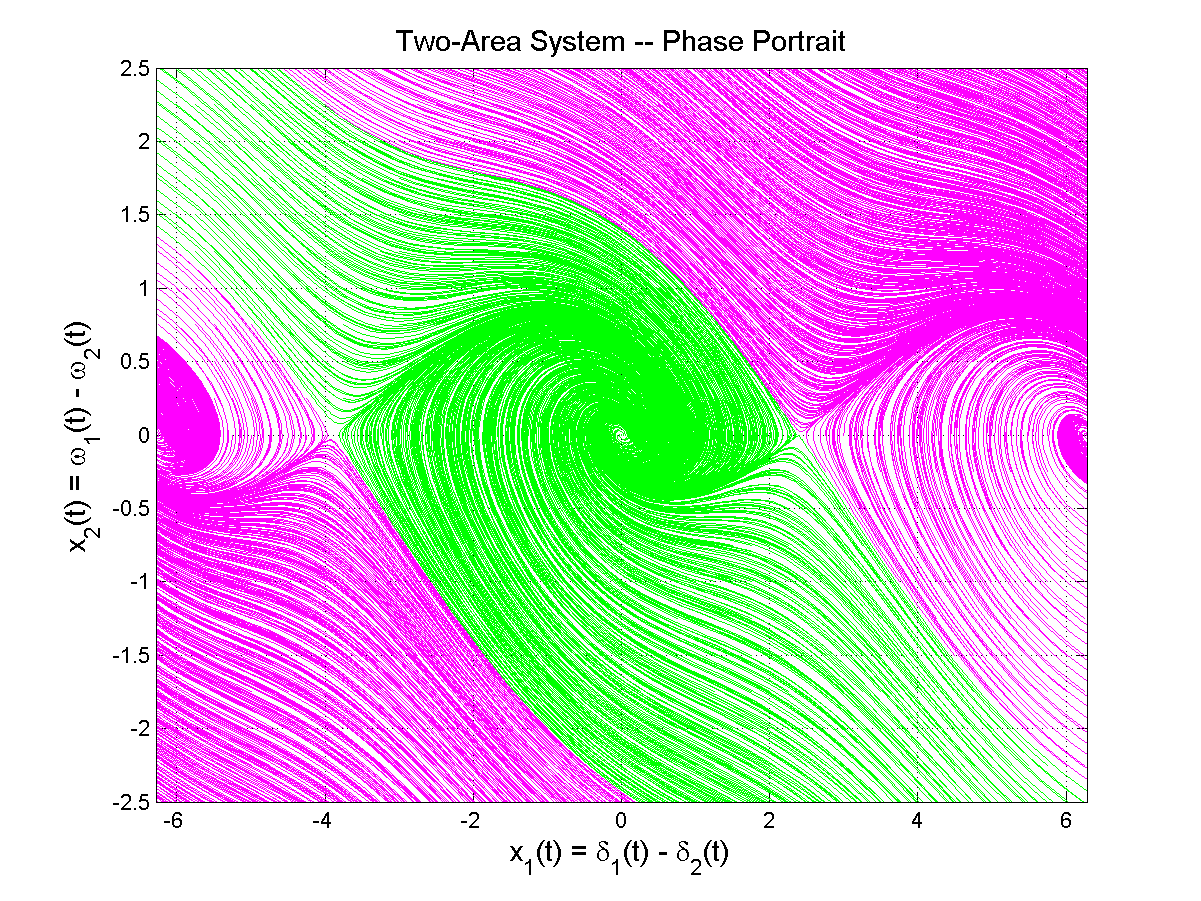}}
{\includegraphics[trim = 1.0cm   0cm   1.5cm   0.0cm, clip=true, angle=0, width=0.36\linewidth, keepaspectratio, draft=false]{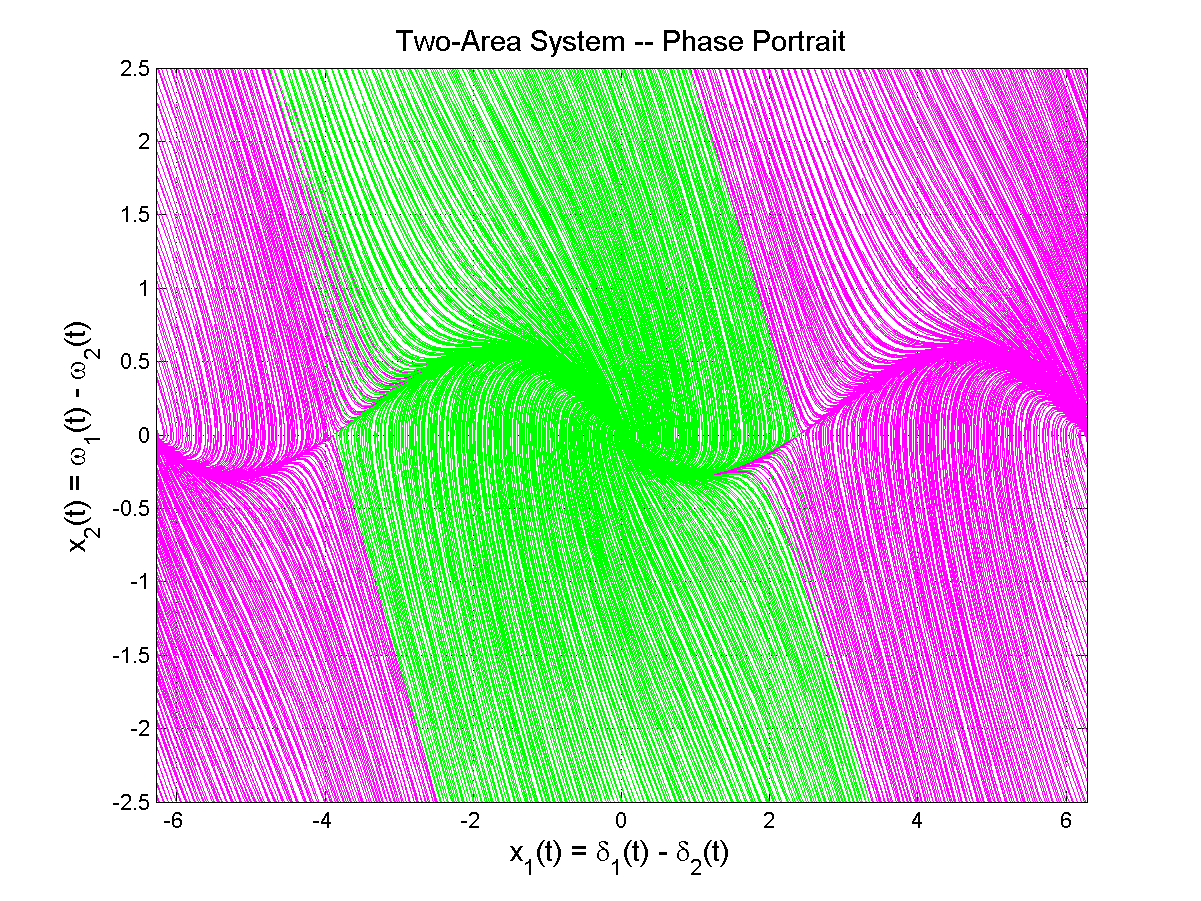}}
{\includegraphics[trim = 1.0cm   0cm   1.5cm   0.0cm, clip=true, angle=0, width=0.36\linewidth, keepaspectratio, draft=false]{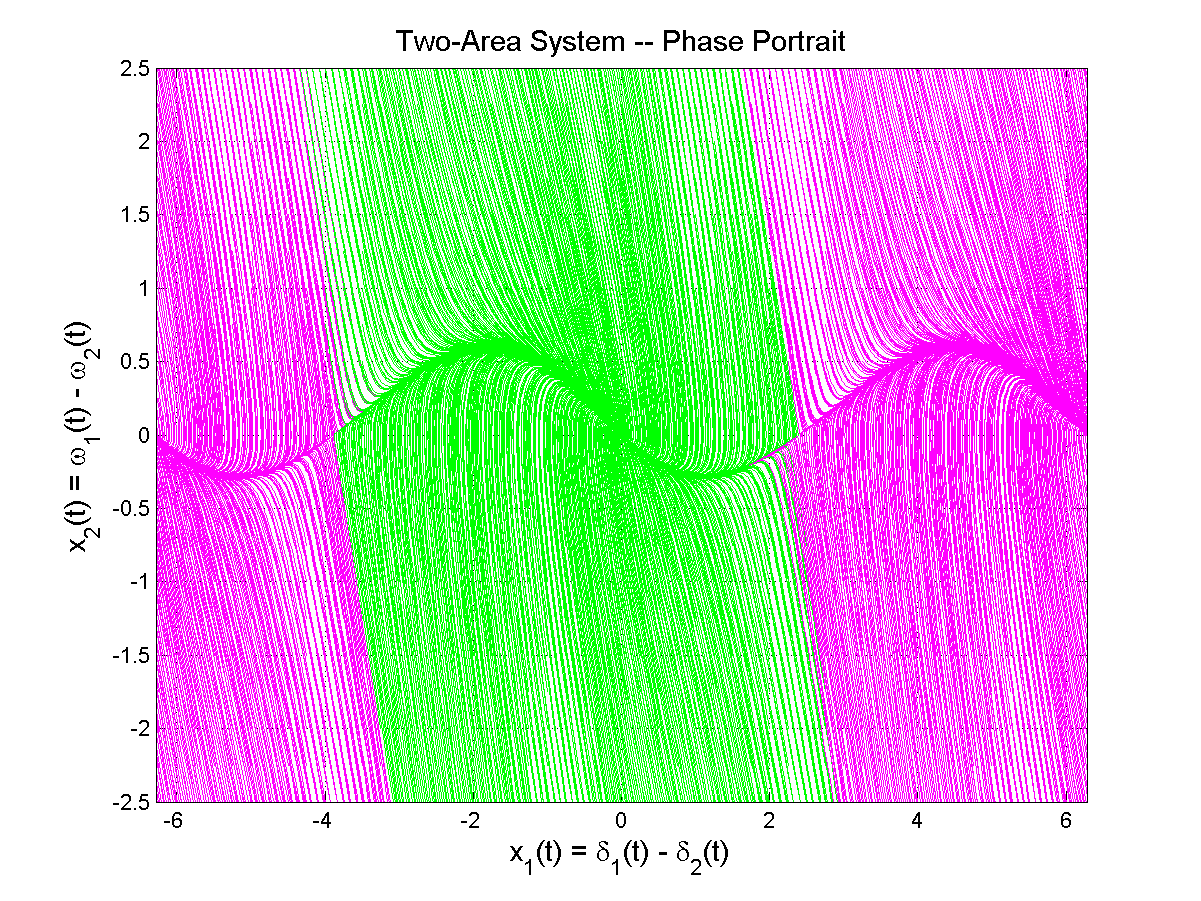}}
  \caption{Unbounded Stability Region of Two-Area System for Different Inertia $M_i$ and Damping $k_i$ (clock-wise). \newline (a)~$M_i = M_0, \ k_i = k_0$, (b)~$M_i = 2 \cdot M_0, \ k_i = k_0$, 
       \newline (c)~$M_i = 0.5 \cdot M_0, \ k_i = 2 \cdot k_0$, (d)~$M_i = M_0, \ k_i = 2 \cdot k_0$.}\label{fig:StabilityRegion}
\end{figure}

In~\cite{Varaiya1987}, the stability region $V(x)$ of a simplified two-area power system around the origin was explicitly calculated and shown. We show in the following that the size and form of this stability region are directly shaped by the choice of the terms $M_i$ and $k_i$. They determine how well \emph{shocks} are absorbed by a power system and how close they drive the system towards the stability boundary $\partial V(x)$. We have calculated the stability region of the two-area system given by Eq.~\ref{eq:Phase} for different choices of $M_i$ and $k_i$. The results are shown in Fig.~\ref{fig:StabilityRegion}. As was stated by~\cite{Varaiya1987}, the stability region (shown in green) is unbounded and centered at the origin.

\begin{figure*}[t]
  \centering               
  \includegraphics{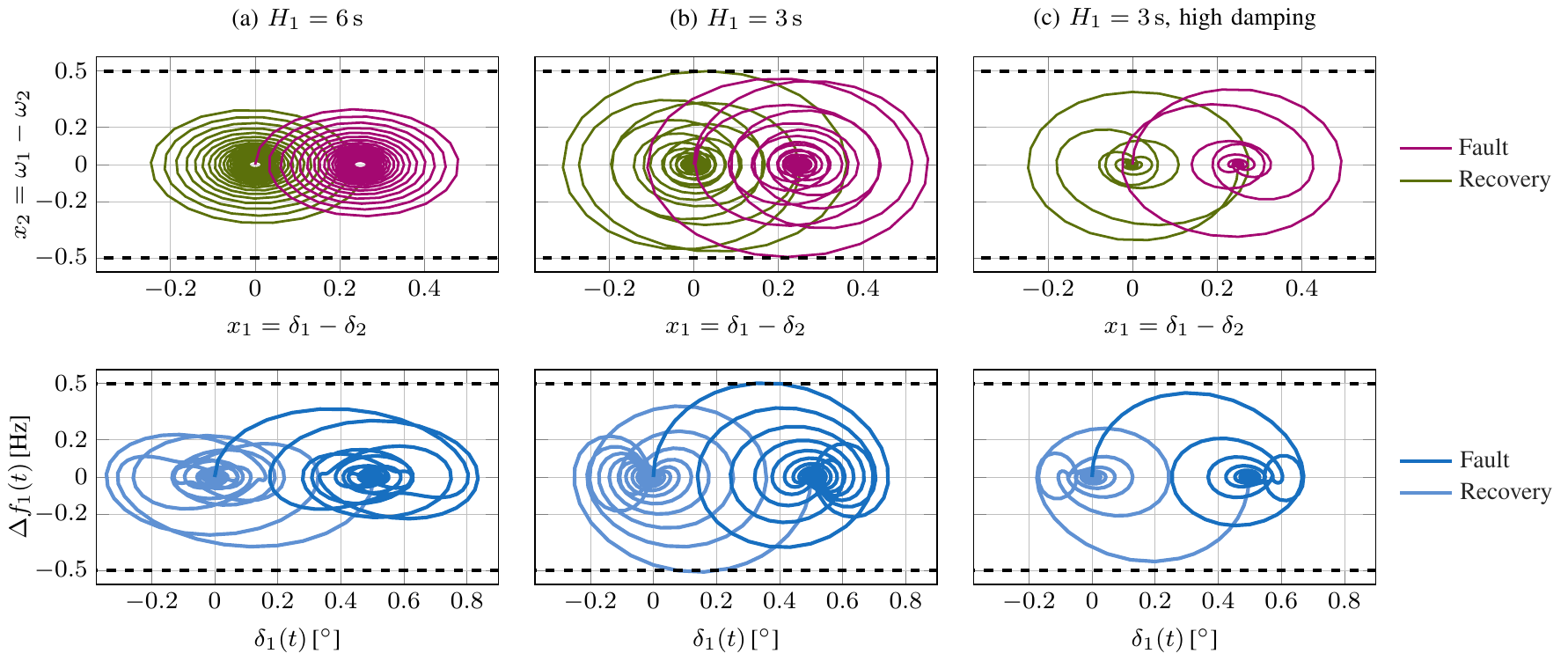}
  \caption{\emph{Upper Plots:} Phase-Plot of Two-Area System, \emph{Lower Plots:} Phase-Plot of Grid Area I. \newline(a)~High Inertia and Low Damping in Grid Area I \footnotesize{($H_1 = H_2 = 6\, \textrm{s}$, $k_1 = k_2 = 1.5\,\frac{\%}{\%}$)}\normalsize.\newline (b)~Low Inertia and Low Damping in Grid Area I \footnotesize{($H_1 = 3\,\textrm{s}$, $H_2 = 6\,\textrm{s}$, $k_1 = k_2 = 1.5\,\frac{\%}{\%}$)}\normalsize. \newline (c)~Low Inertia and High Damping in Grid Area I \footnotesize{($H_1 = 3\,\textrm{s}$, $H_2 = 6\,\textrm{s}$, $k_1 = 4.5\,\frac{\%}{\%}$, $k_2 = 1.5\,\frac{\%}{\%}$)}\normalsize.}\label{fig:PhasePlots}
\end{figure*}

The stability region extends along two axes: the frequency angle difference $x_1 := \delta_1 - \delta_2$ and the frequency deviation $x_2 := f_1 - f_2 = \frac{1}{2 \pi}\left(\omega_1 - \omega_2\right)$.
It is of critical importance that this region is sufficiently large along the $x_2$-axis, since any power fault event happening in a grid region $i$ itself ($\Delta P_i$) or imported via the power lines from neighboring grid regions $j$ ($\Delta P_{i,j}^{\textrm{tie}}$) has a direct impact on $\omega_i$. Rotational inertia is beneficial in reducing the direct impact of $\Delta\,P_i$ on $\omega_i$, i.e.~the excursion of the system state from the origin along the $x_1$-axis, whereas the damping coefficient $k_i$ is good for increasing the size of the stability region along the $x_2$-axis. As we will show in the next section, additional damping can be emulated by fast primary frequency control.

Illustrations of the effect of different values of $M_i$ and $k_i$ are given for the two-area power system in the form of phase-plots (Fig.~\ref{fig:PhasePlots} -- upper plots). Here, the impact of a \emph{shock}, i.e.~a power deviation in Grid Area I given by $\Delta P_i$, is simulated. This results in an excursion of the system state away from the origin to a new equilibrium point on the $x_1$-axis, as (trajectory shown in magenta). After a while, the power fault is cleared and the system then moves back towards the origin (trajectory shown in green). Depending on the choice of parameters $M_i$ and $k_i$ the critical excursion of the system's phase trajectory along the $x_2$-axis is smaller (for large values of $M_i$ and $k_i$) or larger (for small values of $M_i$ and $k_i$). Note that frequency deviations $\Delta f$ of more than $\pm 0.5\textrm{Hz}$ may cause considerable generation tripping.
The additional phase-plot trajectories of Grid Area I, (~Fig.~\ref{fig:PhasePlots} -- lower plots), show that this critical limit is indeed violated in one instance~(Fig.~\ref{fig:PhasePlots} -- bottom, center).

\section{Impact of Low Rotational Inertia on Power System Operation}\label{sec:Operation}

Besides the more theoretical power system stability analysis of the previous chapter, we have also identified impacts of low rotational inertia on daily operational practices in the power systems domain. 

In power systems in general, faster frequency dynamics due to lower levels of rotational inertia raise the question whether fast frequency control, e.g~the primary frequency control scheme in the continental European grid area of ENTSO-E, will remain sufficiently \emph{fast} for mitigating fault events before a critical frequency drop can occur. 
In interconnected power systems in particular, faster frequency dynamics also mean that the swing dynamics of the individual grid areas with their neighboring grid areas will likely be amplified, which in turn leads to significantly amplified transient power exchanges over the power lines.

In current practice stable power system operation is provided by traditional frequency control, which in~\cite{ENTSOE2009} has three categories: Primary frequency control is provided within a few seconds, usually~30~s, after the occurrence of a frequency deviation. It provides power output proportional to the deviation $\Delta f$ ($u_{\textrm{prim.}} = -\frac{1}{S}\, \Delta f$), stabilizing the system frequency but not restoring it to $f_0$. 
Generators of all grid control zones are participating in primary control. The responsible units in the control zone of the imbalance start to take over after approximately $30\,\textrm{s}$, providing secondary frequency control. As secondary control has an integral control part (PI control), it restores both the grid frequency from its residual deviation and the corresponding tie-line power exchanges with other control zones to the set-point values. Tertiary frequency control manually adapts power generation and load set-points and allows the provision of control reserves for grid operation beyond the initial 15~minute time-frame after a fault event has occurred. In addition, generator and load rescheduling can be  manually activated according to the expected residual fault in order to relieve tertiary control by cheaper sources at a later stage, i.e.~with a delay of 45--75 minutes. 


\subsection{Experiments with a One-Area Power System Model}

Due to the faster frequency dynamics, fault events, i.e.~power deviations, have a higher impact on power systems during low rotational inertia situations than usual~\cite[Fig.~14]{Karlsson2008}. We illustrate this by analyzing the dynamic response of the Continental European area power system to fault events, including the stabilizing effect of primary and secondary frequency control schemes. 

\begin{figure}[tbh]
  \centering               
  {\includegraphics[trim = 0.0cm   0.0cm   0.0cm   0.0cm, clip=true, angle=0, width=0.75\linewidth, keepaspectratio, draft=false]{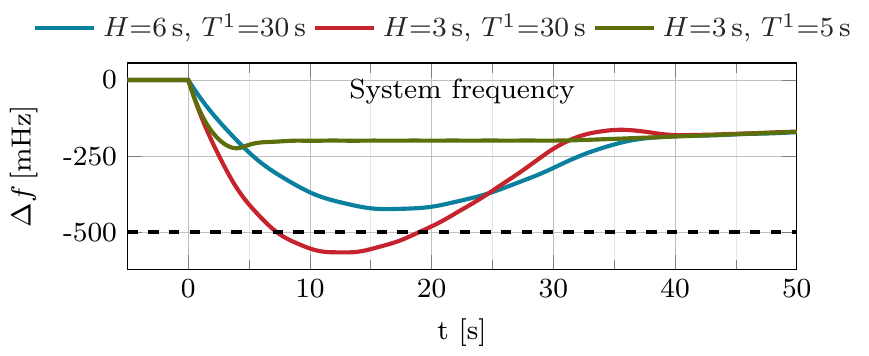}}
  \caption{Dynamic response of the Continental European area power system to faults~\cite{ENTSOE2009}.\vspace{0.15cm}\newline 
  \textcolor{eth3}{\textbf{Blue}}: high inertia ($H=6\,\textrm{s}$), i.e.~no wind\&\ac{pv} power feed-in share, nominal frequency control reserve.\newline
  \textcolor{eth7}{\textbf{Red}}: low inertia ($H=3\,\textrm{s}$), i.e.~$50\,\%$ wind\&\ac{pv} power feed-in share, nominal frequency control reserve.\newline
  \textcolor{eth4}{\textbf{Green}}: low inertia ($H=3\,\textrm{s}$), fast control reserves.\newline
  }
  \label{fig:OneArea}
  \vspace{-0.250cm}
\end{figure}

An Aggregated Swing Equation (ASE), as introduced in Eq.~\ref{eq:ASE}, is considered. Realistic system parameters as identified from actual measurements of the interconnected European system were taken from~\cite{Weissbach2008}. 
A typical summer load demand situation is assumed, e.g.~$230\,\textrm{GW}$ (15 August 2012, 8--9am MEST), and different values of the inertia constant~$H$ are considered. The design worst-case power fault event, an abrupt loss of $\Delta P = 3000\,\textrm{MW}$, is applied to the power system. 
Nominal primary and secondary frequency control schemes are employed, i.e.~primary frequency control reacts with a maximum delay of $5\,\textrm{s}$ and shall achieve full activation after $30\,\textrm{s}$. This corresponds exactly to the control reserve requirements as stated by~\cite{ENTSOE2009}. As shown in Fig.~\ref{fig:OneArea}, the design worst-case power fault event that the continental European system should still be able to sustain, can be absorbed successfully as expected during a high inertia situation ($H_{\textrm{agg}}= 6\,\textrm{s}$) (trajectory shown in black). However, the same fault event becomes critical during a low inertia situation ($H_{\textrm{agg}}= 3\,\textrm{s}$) since the system frequency drops below 49.5~Hz (trajectory shown in red) before the nominal primary frequency control fully kicks in~(30~s~after the fault). In this case the automatic shedding of a combined wind\&\ac{pv} capacity well above 10~GW is, in the current power system setup (year 2013), not merely a theoretical but rather a likely possibility due to the currently existing grid code regulations regarding the fault-ride through behavior of these units.

As can also be seen in this simulation example (shown in green), one powerful mitigation option for low inertia levels and faster frequency dynamics is the deployment of a faster primary control scheme, e.g.~fully activated within 5~s after a fault. Notably Battery Energy Storage Systems (BESS) are well-suited for providing a fast power response as was shown in~\cite{Berlin1986}, \cite{Oudalov2007}, \cite{Ulbig2010c} and~\cite{Borsche2013GM}. Another viable option is the provision of temporary primary frequency control from (variable speed) wind turbines~\cite{Karlsson2008}.
Such a fast primary control response can be thought of as an additional damping term $k_{{\textrm{prim.}}} = \frac{1}{S}$ for the power system as is illustrated by Eq.~(\ref{eq:Damping2}). This effect, depending on its reaction time and power ramping constraints, may provide a crucial stabilization effect in the first seconds after a fault event $\Delta P$. 
This relationship is as follows 
\begin{eqnarray}\label{eq:Damping1}
	\dot{x}	&=& Ax + B_u u + B_d d \, , \ u := -K x \nonumber \\
	\dot{x} &=& Ax + B_u \left(  -K x \right) + B_d d = \left( A - B_u K \right) x + B_d d \nonumber \\
    \Delta\dot{f} &=& A \Delta f + B_u u_{\textrm{prim.}}+ B_d \Delta P \, ,\ \ u_{\textrm{prim.}} := -\frac{1}{S} \quad ,
\end{eqnarray}

where the term $u$ is the control input, i.e.~$u_{\textrm{prim.}} = -\frac{1}{S}$ with $S$ as the bias of the primary frequency control, $d$ a disturbance, i.e.~a power fault event $\Delta P$, and $x = \Delta f$ the system state, i.e.~the grid frequency deviation. 

With $A = -\frac{f_0}{2 H S_\textrm{B}} \cdot k_{\textrm{load}} = -\frac{f_0}{2 H S_\textrm{B}} \cdot \frac{1}{D_\textrm{load}}$, $B_u = B_d = \frac{f_0}{2 H S_\textrm{B}}$ this finally leads to    
\begin{eqnarray} \label{eq:Damping2}
	\Delta\dot{f} &=& \frac{f_0}{2 H S_\textrm{B}} \cdot \left( \underbrace{\left( -\frac{1}{D_\textrm{load}} \right) \Delta f}_{\textrm{Load Damping}} \  + \underbrace{\left( -\frac{1}{S} \Delta f \right)}_{\textrm{Prim. Freq. Ctrl.}} + \ \Delta P \right), \nonumber \\
\Delta\dot{f} &=& \frac{f_0}{2 H S_\textrm{B}} \cdot \left( \underbrace{- \left( k_{\textrm{load}} + k_{\textrm{prim.}}(t) \right) \cdot \Delta f}_{\textrm{Augmented Frequency Damping}} + \underbrace{\Delta P}_{\textrm{Fault}} \right) \, .
\end{eqnarray}

\begin{figure*}[h]
  \centering\includegraphics[width=0.85\textwidth]{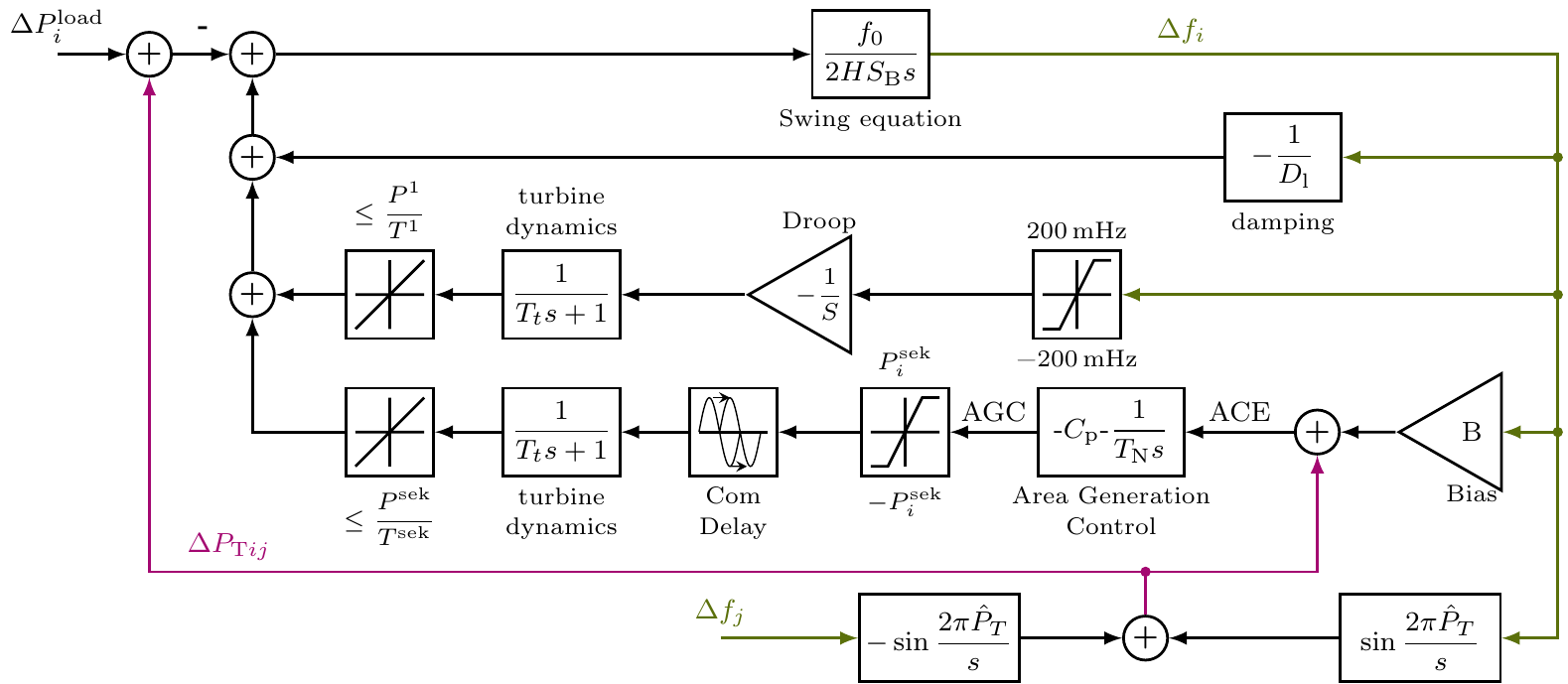}
  \caption{Generalized Multi-Area System (only Grid Area $i$ shown). Implementation in \texttt{Matlab/Simulink}.} 
  \label{fig:Simulink}
\end{figure*}

Note that due to the time-delay behavior of primary frequency control, i.e.~$u_{\textrm{prim.}}(t) = -\frac{1}{S}\, \Delta f(t-T_\textrm{delay})$ and power ramp-rate limitations as shown in Fig.~\ref{fig:Simulink}, the damping effect of the primary frequency control in reality turns out to be a more complex time-variant term, i.e.~$k_{{\textrm{prim.}}}(t)$.

The above swing dynamics (Eq.~\ref{eq:Damping2}) clearly show that the two principal design options for mitigating the impact of power imbalance faults ($\Delta P$) on grid frequency disturbances ($\Delta f$) are to either increase the rotational inertia constant $H$ and/or augment the frequency damping via the provision of fast primary frequency control $k_{{\textrm{prim.}}}(t)$.

\subsection{Experiments with a Two-Area Power System Model}

Unlike to a One-Area system model, which is assumed to represent highly meshed and thus highly coupled grid areas, noticable swing dynamics are observable between more loosely coupled grid areas. An illustration of this is given in the following for a Two-Area power system that shall represent again the continental European power system. The two grid areas are equal in size, their sum being equivalent to the actual system size of the continental European system. We have tried to model the system as realistically as possible, again using the parameters identified in~\cite{Weissbach2008} as well as by incorporating primary and secondary frequency control schemes as illustrated for a generalized, nonlinear multi-area power system in Fig.~\ref{fig:Simulink}. Furthermore, realistic delay, power ramping and saturation blocks are included.

\begin{figure*}[t]
  \centering               
  \includegraphics[trim = 0.0cm  0.0cm 0.0cm 0.0cm, clip=true, angle=0, height=10cm, draft=false]{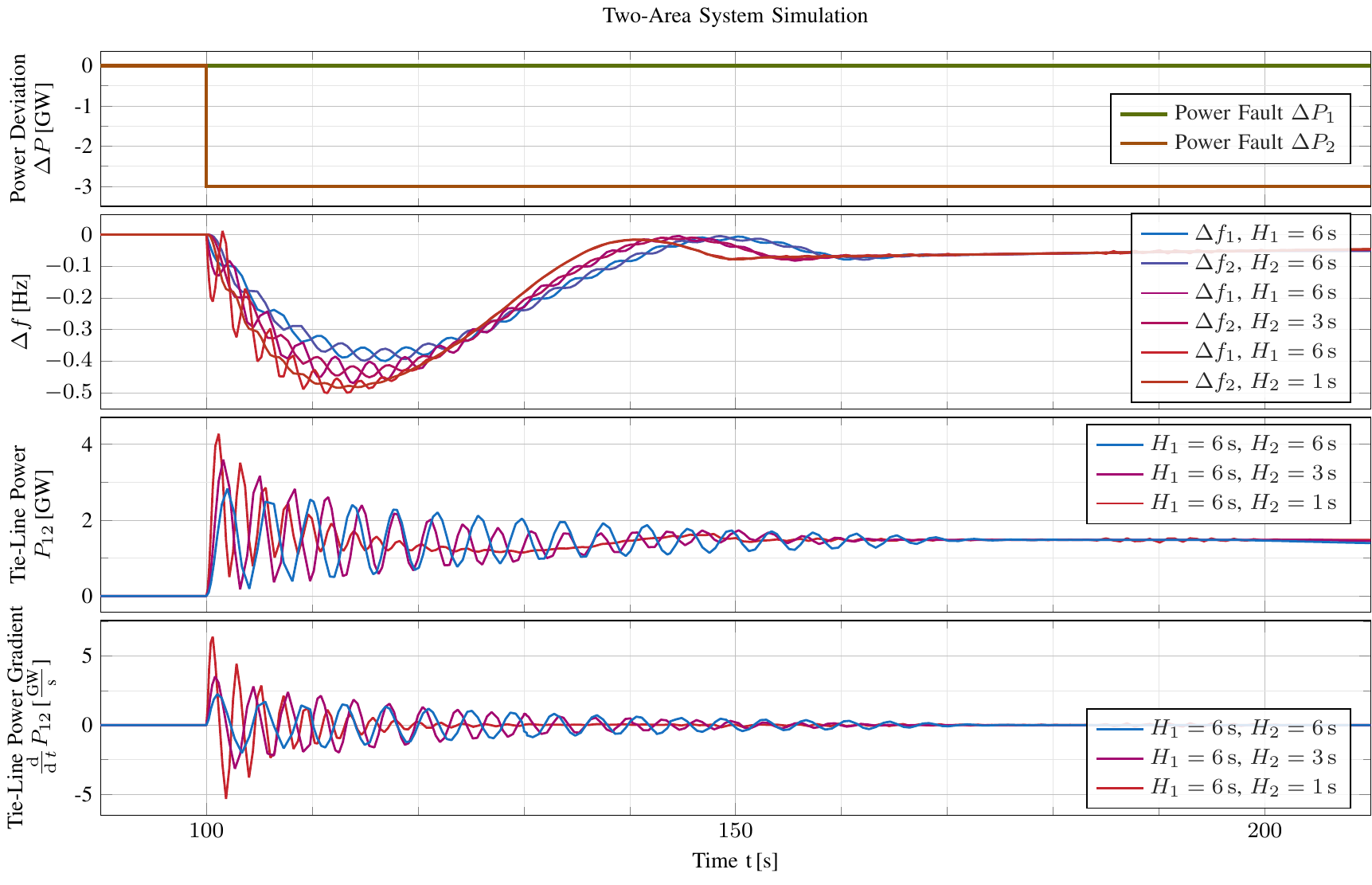}
  \vspace{-0.25cm}
  \caption{Dynamic response of Two-Area System for design worst-case fault (sudden loss of 3000~MW)~\cite{ENTSOE2009}.}  
  \label{fig:TwoArea}
\end{figure*}

In the subsequent simulations, we chose a similar setup as before and again the design worst-case power fault event (3000~MW) occurring after 100~s into the simulation runs. We assumed different levels of rotational inertia in Grid Area II, $H_{\textrm{II}} = \{\, 1\,\textrm{s} \, , \, 3\,\textrm{s} \, , \, 6\,\textrm{s}\, \}$, whereas the rotational inertia in Grid Area I is nominal ($H_{\textrm{I}} = 6\,\textrm{s}$), the base power is split equally ($S_{\textrm{B,I}} = S_{\textrm{B,II}} = 115\,\textrm{GW}$) and everything else is the same.  
The simulation results~(Fig.~\ref{fig:TwoArea}) show that indeed noticeable frequency swing dynamics are observable between the two regions. The swing dynamics are more amplified for lower inertia levels in Grid Area II. As a consequence, the transient power flows $\Delta P_{\textrm{I,II}}^{\textrm{tie}}$ over the tie-line between Grid Areas I and II are significantly increased (by more than 50\%) and becoming more abrupt (by up to 300\%). 
Both the magnitude of transient tie-line power flows as well as their time-derivative $\Delta \dot{P}_{\textrm{I,II}}^{\textrm{tie}}$ can be triggers for automatic protection devices that are designed to clear short circuits by tripping tie-lines. 
In a grid situation as described here, a false short circuit event may be detected by protection devices leading to the immediate tripping of the tie-line in an already sensible moment.

Supplementary experiments with a Three-Area power system show that the phenomenon of swing dynamics and large transient power flows on the tie-lines diminishes, the better meshed the overall system is, i.e.~the more tie-lines exist between the grid areas. Here two possible grid setups exist: connection of the three areas either in the form of a string or a (better meshed) triangle. In the latter case the size of the swing dynamics and transient power flows are smaller and better damped but still remain significant.


\begin{table}[tbh]
  \caption[width=5cm]{Power System Model Parameters.}
  \label{tab:parameters}
  \centering
  \begin{tabular}{l l r@{\,}l r@{\,}l}
    \toprule
    Parameter                & Variable               & \multicolumn{2}{c}{Grid Area I} & \multicolumn{2}{c}{Grid Area II} \\ \midrule
    Rotational Inertia                  & $H$                    & 6    & \si{\second}    & 1$\vert$3$\vert$6     & \si{\second} \\
    Damping                  & $k_{\textnormal{load}}$ & 1.5 & $\frac{\%}{\%}$  & 1.5 & $\frac{\%}{\%}$ \\
    Base Power               & $S_\textnormal{B}$     & 115    & \si{\giga\watt} & 115   & \si{\giga\watt} \\
    Tie-Line Power Rating    & $\hat{P}_\textnormal{T}$     & 0.025   & $S_\textnormal{B}$ & 0.025   & $S_\textnormal{B}$ \\
    Primary Control 		 & $P^\textnormal{1}$  &  1500  & \si{\mega\watt} & 1500  & \si{\mega\watt} \\
    Prim. Response Time    & $T^1$                  & 30   & \si{\second}    & 5$\vert$30    & \si{\second} \\
    Secondary Control 		 & $P^\textrm{sek}$ & 14000  & \si{\mega\watt} & 14000 & \si{\mega\watt} \\
    Sec. Response Time  & $T^\textrm{sek}$                  & 120  & \si{\second}    & 120   & \si{\second} \\
    AGC Parameters           & $C_\textnormal{p}$     & 0.17 &                 & 0.17  &  \\
                             & $T_\textnormal{N}$     & 120  & \si{\second}    & 120   & \si{\second} \\
\bottomrule
  \end{tabular}
\end{table}

\pagebreak
\vspace{-0.05cm}
\section{Conclusion and Outlook}\label{sec:Conclusion}
\vspace{-0.20cm}
The presented analyses show that high shares of inverter-connected power generation can have a significant impact on power system stability and power system operation.\newline
The new contributions of this paper are:
\begin{itemize}
 	\item \textbf{Rotational Inertia becomes heterogeneous.} 
 	Instead of a global inertia constant $H$ there are different $H_i$ for the individual areas $i$ as a function of how much converter-connected units versus conventional units are online in the different areas.
 	\item \textbf{Rotational inertia constants become time-variant ($H_i(t)$).} 
 	This is due to the variability of the power dispatch. Frequency dynamics become thus differently fast in the individual grid areas.
 	\item \textbf{Grid frequency instability phenomena are amplified.} 
 	Reduced rotational inertia leads to faster frequency dynamics and in turn causes larger frequency deviations and transient power exchanges over tie-lines in the event of a power fault. This may cause false errors and unexpected tripping of the tie-lines in question by automatic protection devices, in turn further aggravating an already critical situation.
 	\item \textbf{Faster primary control emulates a time-variant damping effect ($k(t)$).} 
 	This is critical for power system stability immediately after a fault event.
\end{itemize}

Please note that the analysis results presented here have been obtained by using \emph{idealized} primary and secondary frequency control loop dynamics. This is only a first step. Further analysis will, however, have to take into account more detailed, i.e.~more realistic, frequency response characteristics of various unit types (i.e. including additional time-delays, inverse response behavior, etcetera).

Mitigation options for low rotational inertia and faster frequency dynamics are faster primary frequency control and the provision of synthetic rotational inertia, also known as inertia mimicking, provided either by wind\&\ac{pv} generation units and/or storage units; confer also to~\cite{Karlsson2008}, \cite{Borsche2013GM}, \cite{MULLANE_2005}, \cite{MORREN_2006b} and \cite{Ulbig2013CDC}. 
\newline
BESS units are, due to their very fast response behavior, especially well-suited for providing either fast frequency (and voltage) control reserves or synthetic rotational inertia for power system operation.
\newpage

\vspace{-0.2cm}                 
\bibliographystyle{IEEEtran}	
\bibliography{Inertia_FullPaper}

\begin{thebibliography}{10}
\providecommand{\url}[1]{#1}
\csname url@samestyle\endcsname
\providecommand{\newblock}{\relax}
\providecommand{\bibinfo}[2]{#2}
\providecommand{\BIBentrySTDinterwordspacing}{\spaceskip=0pt\relax}
\providecommand{\BIBentryALTinterwordstretchfactor}{4}
\providecommand{\BIBentryALTinterwordspacing}{\spaceskip=\fontdimen2\font plus
\BIBentryALTinterwordstretchfactor\fontdimen3\font minus
  \fontdimen4\font\relax}
\providecommand{\BIBforeignlanguage}[2]{{%
\expandafter\ifx\csname l@#1\endcsname\relax
\typeout{** WARNING: IEEEtran.bst: No hyphenation pattern has been}%
\typeout{** loaded for the language `#1'. Using the pattern for}%
\typeout{** the default language instead.}%
\else
\language=\csname l@#1\endcsname
\fi
#2}}
\providecommand{\BIBdecl}{\relax}
\BIBdecl

\bibitem{Kundur1994}
P.~Kundur, ``Power system stability and control,'' \emph{McGraw-Hill Inc., New
  York}, 1994.

\bibitem{REN21:2012}
\BIBentryALTinterwordspacing
REN21, ``Renewables 2012 global status report,'' 2012. [Online]. Available:
  \url{www.ren21.net}
\BIBentrySTDinterwordspacing

\bibitem{RESfigure}
\BIBentryALTinterwordspacing
BMU, ``{Renewable Energy Sources in Figures -- National and International
  Development},'' 2013. [Online]. Available:
  \url{www.erneuerbare-energien.de/files/english/pdf/application/pdf/broschuere_ee_zahlen_en_bf.pdf}
\BIBentrySTDinterwordspacing

\bibitem{VanHertem2012}
P.~Tielens and D.~Van~Hertem, ``Grid inertia and frequency control in power
  systems with high penetration of renewables,'' \emph{Young Researchers
  Symposium in Electrical Power Engineering, Delft}, vol.~6, April 2012.

\bibitem{Kopell1982}
N.~Kopell and R.~B.~J. Washburn, ``Chaotic motions in the two-degree-of-freedom
  swing equations,'' \emph{IEEE Transactions on Circuits and Systems}, vol.~29,
  no.~11, pp. 738--746, 1982.

\bibitem{Varaiya1987}
H.-D. Chiang, F.~F. Wu, and P.~P. Varaiya, ``{Foundations of Direct Methods for
  Power System Transient Stability Analysis},'' \emph{IEEE Transactions on
  Circuits and Systems}, February 1987.

\bibitem{Berggren1993}
B.~Berggren and G.~Andersson, ``On the nature of unstable equilibrium points in
  power systems,'' \emph{Power Systems, IEEE Transactions on}, vol.~8, no.~2,
  pp. 738--745, 1993.

\bibitem{ENTSOE2009}
\BIBentryALTinterwordspacing
ENTSO-E, ``Operation {H}andbook,'' 2009. [Online]. Available:
  \url{www.entsoe.eu/resources/publications/entso-e/operation-handbook/}
\BIBentrySTDinterwordspacing

\bibitem{Karlsson2008}
N.~Ullah, T.~Thiringer, and D.~Karlsson, ``{Temporary Primary Frequency Control
  Support by Variable Speed Wind Turbines -- Potential and Applications},''
  \emph{Power Systems, IEEE Transactions on}, vol.~23, no.~2, pp. 601--612,
  2008.

\bibitem{Weissbach2008}
T.~Weissbach and E.~Welfonder, ``{Improvement of the Performance of Scheduled
  Stepwise Power Programme Changes within the European Power System},'' in
  \emph{17th {IFAC} {W}orld {C}ongress, {T}he International Federation of
  Automatic Control ({IFAC})}, Seoul, Korea, 2008, pp. 11\,972--11\,977.

\bibitem{Berlin1986}
H.~Kunisch, K.~Kramer, and H.~Dominik, ``Battery energy storage another option
  for load-frequency-control and instantaneous reserve,'' \emph{Energy
  Conversion, IEEE Transactions on}, vol. EC-1, no.~3, pp. 41--46, 1986.

\bibitem{Oudalov2007}
A.~Oudalov, D.~Chartouni, and C.~Ohler, ``{Optimizing a Battery Energy Storage
  System for Primary Frequency Control},'' \emph{IEEE Transactions on Power
  Systems}, vol.~22, no.~3, pp. 1259--1266, Aug. 2007.

\bibitem{Ulbig2010c}
A.~Ulbig, M.~D. Galus, S.~Chatzivasileiadis, and G.~Andersson, ``General
  {F}requency {C}ontrol with {A}ggregated {C}ontrol {R}eserve {C}apacity from
  {T}ime-{V}arying {S}ources: {T}he {C}ase of {PHEV}s,'' in \emph{IREP
  Symposium 2010 -- Bulk Power System Dynamics and Control -- VIII}, Buzios,
  RJ, Brazil, 2010.

\bibitem{Borsche2013GM}
T.~Borsche, A.~Ulbig, M.~Koller, and G.~Andersson, ``{Power and Energy Capacity
  Requirements of Storages Providing Frequency Control Reserves},'' in
  \emph{IEEE PES General Meeting}, Vancouver, 2013.

\bibitem{MULLANE_2005}
A.~Mullane and M.~O'Malley, ``The inertial response of induction-machine-based
  wind turbines,'' \emph{IEEE Transactions on Power Systems}, vol.~20, no.~3,
  Aug. 2005.

\bibitem{MORREN_2006b}
J.~Morren, S.~de~Haan, and J.~Ferreira, ``Contribution of {D}{G} units to
  primary frequency control,'' in \emph{2005 {I}nternational {C}onference on
  {F}uture {P}ower {S}ystems}, Nov. 2005.

\bibitem{Ulbig2013CDC}
A.~Ulbig, T.~Rinke, S.~Chatzivasileiadis, and G.~Andersson, ``{Predictive
  Control for Real-Time Frequency Regulation and Rotational Inertia Provision
  in Power Systems},'' in \emph{accepted at CDC 2013}, Florence, 2013.

\end{thebibliography}

\end{document}